\documentclass[12pt]{article}
\usepackage{amssymb}
\usepackage{graphicx}
\usepackage{latexsym}
\def\part#1{\frac{\partial\phantom{q}}{\partial#1}}

\newenvironment{rmk}{\begin{trivlist}\item[]{\bf Remark:} }
{\end{trivlist}}

\newenvironment{prf}{\begin{trivlist}\item[]{\bf Proof:} }
{\hfill $\Box$ \end{trivlist}}
\newenvironment{lemprf}{\begin{trivlist}\item[]{\bf Proof:} }
 {\end{trivlist}}
\newtheorem{thm}{Theorem}

\newtheorem{prp}[thm]{Proposition}
\newtheorem{lemma}[thm]{Lemma}

\newcommand{\lie}[1]{\mathfrak{#1}}

\def\deg{\mathop{\rm deg}\nolimits}

\def\tr{\mathop{\rm tr}\nolimits}
\def\td{\mathop{\rm td}\nolimits}

\def\ch{\mathrm {ch}}

\newcommand{\R}{\mathbf{R}}
\newcommand{\C}{\mathbf{C}}

\newcommand{\Z}{\mathbf{Z}}

\newcommand{\PP}{{\rm P}}

\begin{document}
\title{Spherical harmonics and the icosahedron}
 \author{Nigel Hitchin\\[5pt]}
\maketitle
\centerline{\it Dedicated to John McKay}
\section{Introduction}
Spherical harmonics of degree $\ell$  are the functions on the unit $2$-sphere which satisfy $\Delta f=\ell(\ell+1)f$ for the Laplace-Beltrami operator $\Delta$. They form an irreducible representation of $SO(3)$ of dimension $2\ell+1$ and are the restrictions to the sphere of homogeneous polynomials $f(x_1,x_2,x_3)$ of degree $\ell$ which solve Laplace's equation in $\R^3$. This paper concerns a curious relationship between the case $\ell=3$ and the regular icosahedron. 

We consider the zero set of $f$ -- the {\it nodal set} -- and ask whether it contains the $12$ vertices of a regular icosahedron. If, for example,  $f=(5 \cos^3 \theta-3\cos \theta)/2$ (the third order Legendre polynomial) then its nodal set consists of the intersection of the sphere with the three planes  $x_3=0,x_3=\pm \sqrt{3/5}$ and it is easy to see that no three parallel planes can contain all the vertices of the icosahedron. On the other hand, if $\phi=(1+\sqrt {5})/2$, consider the standard icosahedron with   vertices
$(0,\pm \phi, \pm 1)$, $(\pm 1,0,\pm \phi)$, $(\pm \phi,\pm 1,0)$ (on the sphere of radius $\sqrt{\phi^2+1}$). They clearly lie on the nodal set of $f=x_1x_2x_3$, but then so do   
$(0,\pm 1, \pm \phi)$, $(\pm \phi,0,\pm 1)$, $(\pm 1,\pm \phi,0)$ so we have two icosahedra.

We define an invariant of $f$ as follows: let
$$M_{ij}(f)=\int_{S^2}\!\!f^2(x)(x_ix_j-\frac{4}{15}\delta_{ij})\,dS$$
and let $x^3+a_1 x^2+a_2x+a_3$ be the characteristic polynomial of $M$. Put
$J(f)=3a_3-4a_1a_2$, then our  result is:
\newpage
{\bf Theorem}: {\it Let $f$ be a degree three  spherical harmonic.}
\begin{itemize}
\item {\it If $J(f)>0$ (resp.  $J(f)<0$) then the nodal set of $f$ contains the vertices of two (resp. zero) regular icosahedra.} 
\item
{\it When $J(f)=0$ and $f$ is of the form  $f(x)=(u,x)((a,x)^2-(x,x)(a,a)/5)$
with $(u,a)=0$ then any regular icosahedron with $a$ as a vertex lies on the nodal set, otherwise there is a unique such icosahedron.}
\end{itemize}

Our proof uses the geometry of the Clebsch diagonal cubic surface, vector bundles on an elliptic curve and a Fano threefold introduced by S. Mukai.

\section{Harmonic cubics}

Let $U$ be a $3$-dimensional real vector space with positive definite inner product $(a,b)$ and consider the $10$-dimensional vector space  $Sym^3U^*$ of  degree $3$ homogeneous polynomial functions  on $U$. The inner product identifies $U\cong U^*$ and induces an inner product on $Sym^3U^*$ which we   
normalize   so that
$$(p(x),(a,x)^3)=p(a).$$
The special orthogonal group $SO(3)$ acts on $Sym^3U^*$ and decomposes it  into two orthogonal irreducible components: 
a $7$-dimensional  representation $H$ of functions which  satisfy Laplace's equation (the {\it harmonic cubics}), and a $3$-dimensional representation consisting of cubics of the form $(a,x)(x,x)$ for $a\in U$. It is the space $H$ which will concern us here. Restricting to the unit sphere these are eigenfunctions of the Laplace-Beltrami operator with eigenvalue $12$.

If $p\in H$, then the inner product $p(a)=(p(x),(a,x)^3)$ only depends on the harmonic component of $(a,x)^3$.
This is of the form 
$$f(x)=(a,x)^3-\lambda (a,x)(x,x)$$
and differentiating we find that it satisfies Laplace's equation only if $\lambda=3/5$.
Thus 
\begin{equation}
f_a(x)=(a,x)^3-\frac{3}{5} (a,a)(a,x)(x,x)
\label{harm}
\end{equation}
is the harmonic part of $(a,x)^3$ and so 
$$(p(x),f_a(x))=p(a).$$
The function $f_a$ is, up to a constant multiple,  the unique harmonic cubic which is symmetric about the axis $a$ -- if $a=(0,0,1)$ then $f_a=r^3P_3(\cos \theta)$ for the Legendre polynomial $P_3$.

\section{The icosahedron}
Let $\phi=(1+\sqrt{5})/2$ be the  golden ratio. The standard model for a regular icosahedron of side length $2$ has vertices $\pm a_1,\dots,\pm a_6$ given by 
\begin{equation}
(0,\pm \phi, \pm 1),\quad (\pm 1,0,\pm \phi),\quad (\pm \phi,\pm 1,0).
\label{model}
\end{equation}
It has $20$ triangular faces, $30$ edges and $12$ vertices occurring in six opposite pairs. 
\begin{figure}
\begin{center}
\includegraphics [scale=.6] {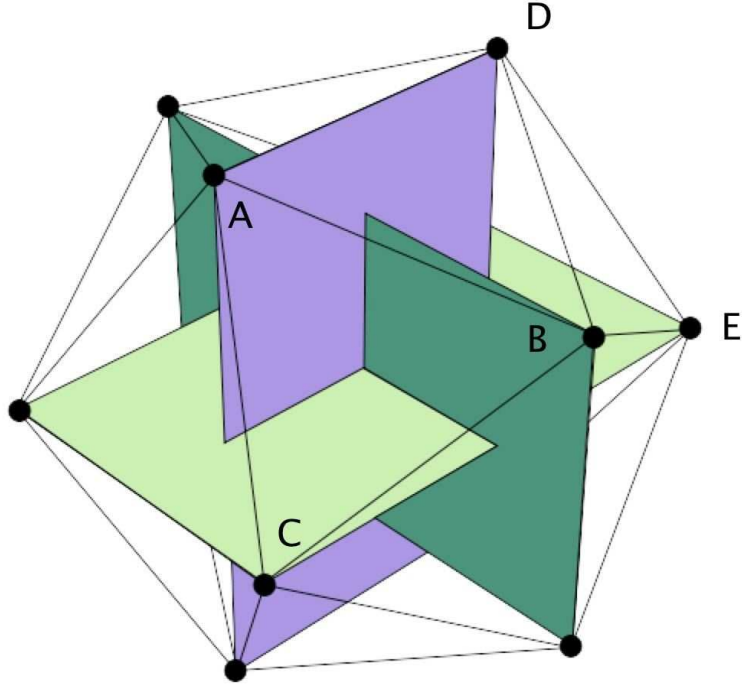}
\end{center}
\caption{ The icosahedron \cite{Fig}}
\end{figure}
The group  of symmetries of the icosahedron is isomorphic to the alternating group $A_5$. The five objects it permutes are triples of orthogonal planes which pass through the vertices (see  Fig.1). In the standard model such a triple is provided by the coordinate planes $x_1=0,x_2=0,x_3=0$ and clearly  the zero set of the cubic $q(x)=x_1x_2x_3$ contains all the vertices of the standard icosahedron.  

  The first result provides the link between the two themes in the title. 

\begin{prp} \label{lem}
The cubic polynomials which vanish on the vertices of a regular icosahedron form a four-dimensional space of  harmonic cubics.
\end{prp}
\begin {prf}  To prove that there is a $4$-dimensional space of such cubic functions  we adopt a  geometric point of view, though a character argument for the icosahedral group can also be used. We use algebraic geometric results  over the complex numbers and then specialize to the real situation. We shall still denote the vector space over the complex numbers by $U$, but recognize that it has a real structure, an antilinear involution.  Then a homogeneous cubic polynomial $p$  defines a cubic curve in the two-dimensional projective plane $\PP(U)$.  The image of the nodal set under the double covering $S^2\rightarrow \PP(U)^r$ consists of the real points of the cubic curve. (We shall use the superscript $r$ to denote real points). 

The cubic curves passing through $m$ points in the projective plane define the anticanonical system for the algebraic surface $S$ obtained by blowing up these points. Take  six points to be defined by the six axes of a regular icosahedron.  Clearly no three are collinear. Also,  no conic passes  through all six points because if it did, its transform by an element of $A_5$ would be another conic meeting it in at least $6$ points and hence coinciding with it. The  conic would then be invariant, and since  $U$ is an irreducible representation of $A_5$  this would be the null conic $(x,x)=0$. But the axis of an icosahedron is not a null vector. 

These facts mean that the anticanonical bundle $K^*$ of $S$ is ample 
and has a four-dimensional space $V$ of sections. Furthermore, since $c_1^2(S)=9-6=3$  these sections  embed $S$ in $\PP^3$ as a nonsingular cubic surface.  

Now if the icosahedron is the standard model, then $q_1(x)=x_1x_2x_3$ defines a   cubic which passes through the six points, as do its  transforms $q_1,\dots,q_5$ by the action of the group $A_5$. These generate the four-dimensional irreducible permutation representation, hence  the $q_{\alpha}$ span this space. The single relation  is $\sum_{\alpha} q_{\alpha}=0$. But $q_1$ and hence $q_{\alpha}$ satisfies Laplace's equation, so all linear combinations of such cubics are harmonic.
\end{prf}

\section{The Clebsch cubic surface}\label{clebsch}

In the course of Proposition \ref{lem} we encountered the cubic surface $S$ obtained by blowing up $\PP(U)$ at the six axes of an icosahedron. The rational map  $\rho:\PP(U)\rightarrow \PP^3$ is described by $x\mapsto [q_1(x),\dots,q_5(x)]$ where $\sum_{\alpha} q_{\alpha}(x)=0$. In fact $\sum_{\alpha} q^3_{\alpha}(x)=0$ and this is the equation of the cubic surface. Indeed, it is invariant in the permutation representation of $A_5$ and since the elementary symmetric function $\sigma_1$ vanishes, there is a unique cubic  invariant polynomial $\sigma_3$, or equivalently the sum of cubes above. This defines  the {\it Clebsch diagonal cubic surface}. Figure 2 shows the plaster model acquired by J.J.Sylvester in 1886 for Oxford University (see \cite{800}). 

\begin{figure}
\begin{center}
\includegraphics [scale=.8] {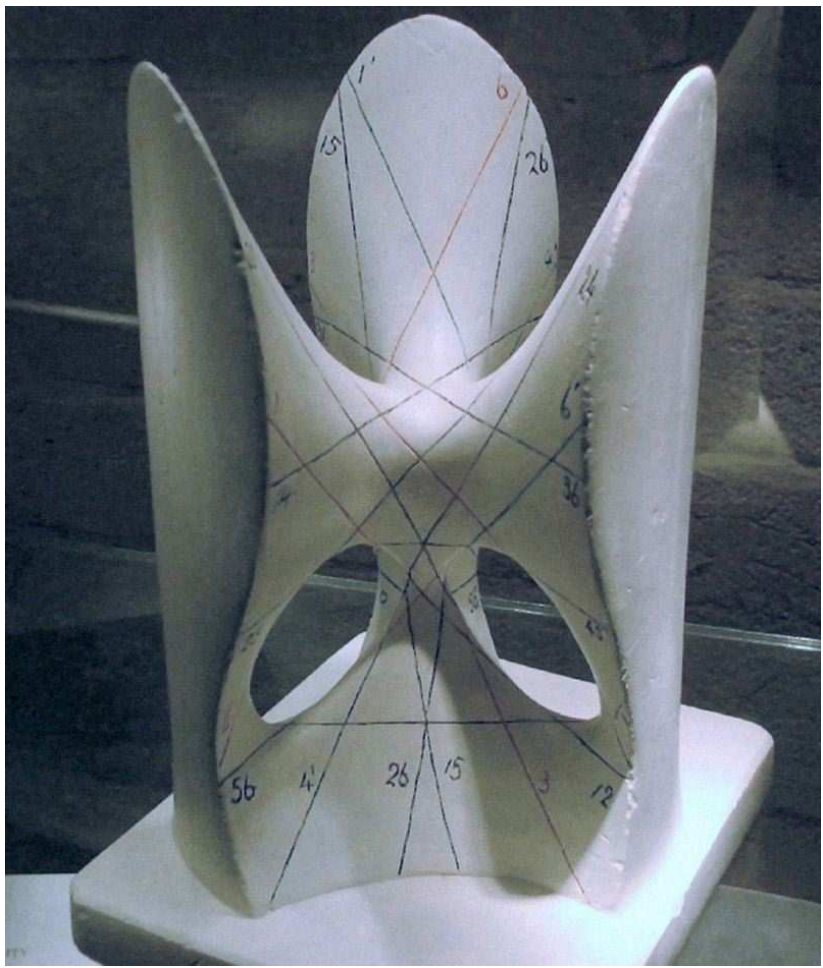}
\end{center}
\caption{ The Clebsch diagonal cubic surface }
\end{figure}

The cubic surface we have just described lies in the projective space $\PP(V^*)$ where $V=H^0(S,K^*)$: given a point $u\in U$ the polynomials in $V$ which vanish at $u$ define a plane in $\PP(V)$. Dually, $[p]\in \PP(V)$ defines a plane in $\PP(V^*)$ and its intersection with the surface $S\subset \PP(V^*)$ blows down to the cubic curve $p=0$ in $\PP(U)$.

However, $V$ has an $A_5$-invariant inner product which identifies it with its dual and so we may regard the Clebsch cubic as also lying in $\PP(V)$. If $p(x)=\sum_{\alpha} y_{\alpha}q_{\alpha}(x)$ with $\sum_{\alpha} y_{\alpha}=0$ then its equation is 
$$\sum_1^5y_{\alpha}^3=0.$$
The rational map  $\rho:\PP(U)\rightarrow S\subset \PP(V)$  has a direct interpretation. Given $u\in U$ let $p_u$ be the orthogonal projection of $f_u$ (defined in (\ref{harm})) onto $V$.  If this projection is zero  then $(p,f_u)=p(u)=0$ for all $p\in V$, so that any cubic which vanishes at the $a_{i}$ also vanishes at $u$. However, these cubics define a projective embedding of the blown-up plane, and so must separate points. It follows that  $p_u=0$ if and only if $[u]=[a_{i}]$ for some $1\le i\le 6$. 

 Now a given $p$ is orthogonal to $f_u$ if and only if  $p(u)=0$, so $p_u\in V$ is orthogonal to the codimension one subspace of $V$ consisting of cubics that vanish at $u$. This, however,  is the definition of the Clebsch surface  in $\PP(V)$. Thus $u\mapsto [p_u]$ maps $\PP(U)\setminus\{[a_1],\dots,[a_6]\}$ to $S$.

Since $p_u$ is the orthogonal projection of $f_u$ on $V$, it follows that 
$(p_u,q_{\alpha})=(f_u,q_{\alpha})=q_{\alpha}(u)$
from which we have an explicit form for the map:
$$\rho(u)=\sum_{\alpha=1}^5 q_{\alpha}(u)q_{\alpha}.$$

There remains the identification within $\PP(V)$ of the six exceptional curves obtained by blowing up the points $[a_{i}]\in \PP(U)$, and the extension of the map $u\mapsto[p_u]$. Consider $p_{u(t)}$ for  $u(t)=a_1+tv$ as $t\rightarrow 0$, where $(v,a_1)=0$.  We have 
$$f_{u(t)}=f_{a_1}+3t(v,x)\left((a_1,x)^2-\frac{1}{5}(a_1,a_1)(x,x)\right)+\dots$$
Consider the  cubic polynomial 
\begin{equation}
h(x)=(v,x)\left((a_1,x)^2-\frac{1}{5}(a_1,a_1)(x,x)\right)
\label{fund}
\end{equation}
where $(v,a_1)=0$. We shall meet this type of harmonic cubic many times. The angle $\theta$  between two different axes of the icosahedron is given by 
$$\cos\theta=\pm \frac{\phi}{1+\phi^2}=\pm \frac{1}{\sqrt{5}}$$
so $h(a_i)=0$ if $i\ne 1$, but also $h(a_1)=0$ since  $(v,a_1)=0$. Thus $q\in V$.  
As we saw above, $p_{a_1}=0$ so $f_{a_1}$ is orthogonal to $V$. Hence $p_{u(t)}=3th+\dots$ and  $v\mapsto h$  defines a map from the line $\PP(a_1^{\perp})$, the projectivized tangent space to $\PP(U)$ at $[a_1]$, to $\PP(V)$. This is the   exceptional curve $E_1$ obtained  by blowing up $\PP(U)$ at $[a_1]$, and is a line in the cubic surface $S\subset \PP(V)$. Geometrically the cubic curve 
$h(x)=0$  is the union of a line through $[a_1]$ and the unique conic $C_1$ passing through the five points $[a_2],\dots,[a_6]$, so the line $E_1$ in $\PP(V)$ consists of the pencil of lines through $[a_1]$ together with $C_1$.

For the six axes we get six disjoint lines --  six of the  $27$ real lines shown on the model in Fig 2. As is well known, the other $21$ lines in the cubic surface are the six proper transforms $E'_1,\dots,E'_6$ of $C_1,\dots,C_6$ together with  the transforms $E_{ij}$ of the $15$ lines $L_{ij}$ joining $[a_i]$ to $[a_j]$. 

The line $E'_i$ consists of the pencil of cubics in $\PP(V)$ with a singularity at $[a_i]$. The   polynomials  $q_{\alpha}-q_{\beta}$,  $1\le \alpha< \beta\le 5$, clearly satisfy $\sum y_{\alpha}=0=\sum y^3_{\alpha}$ and so lie on the Clebsch cubic. There are $15$ disjoint pairs $\{\alpha,\beta\},\{\gamma,\delta\}$ and, as $a,b$ vary, 
$$a(q_{\alpha}-q_{\beta})+b(q_{\gamma}-q_{\delta})$$
describes  $15$ lines $E_{ij}$ in $S$. 

Each $q_{\alpha}-q_{\beta}$ lies in three of these. On a general cubic surface points where three lines intersect are  called Eckard points -- the Clebsch cubic has the maximal number ten of these.

\begin{rmk} Under the composite map $S^2\rightarrow \PP(U)^r\stackrel{\rho}\rightarrow S$, the real lines (apart from the blow-ups $E_i$) are recognizable in the geometry of the icosahedron. For a vertex $a_1$, the adjacent vertices $a_2,\dots,a_6$ lie on a circle in $S^2$ and this maps to $E_1'\subset S$. The points $a_1,a_2,-a_1,-a_2$ lie on a great circle and this maps to  $E_{12}$. The three great circles $E_{12},E_{34},E_{56}$ intersect in the midpoints of a pair of opposite triangular faces, which define the Eckard points.
\end{rmk}

\section{Skew forms}\label{skew}
We now adopt a new viewpoint on our space $H$ of harmonic cubics. 
If $U$ is identified with the Lie algebra of $SO(3)$, the action of a vector $u\in U$ on a homogeneous polynomial is via the differential operator
$$u\cdot f=\sum_{i,j,k}\epsilon_{ijk}x_iu_j\frac{\partial f}{\partial x_k}$$
This is skew-adjoint with respect to the inner product so each $u\in U$ defines a skew form
$$\omega_u(p,q)=(u\cdot p,q)$$
on the $7$-dimensional space $H$. 

Since the Legendre polynomial $f_u$ is symmetric about the axis $u$, clearly $u\cdot f_u=0$. More generally we calculate (using the vector product $a\times b$)
$$u\cdot f_a(x)=3(u\times a,x)\left((a,x)^2-\frac{1}{5}(a,a)(x,x)\right)$$
which is precisely the harmonic cubic (\ref{fund}) we met before.

Any skew form on an odd-dimensional space is degenerate. In the case of $\omega_u$ we have:
\begin{prp} \label{degen}The degeneracy subspace of $\omega_u$ on the space of harmonic cubics is spanned by $f_u$.
\end{prp}
\begin{prf} The degeneracy subspace is defined as the space of cubics  $p\in H$ such that $\omega_u(p,q)=0$ for all $q\in H$. Since
$$\omega_u(p,q)=(u\cdot p ,q)$$
this holds if and only if $u\cdot p=0$, i.e. if $p$ is annihilated by the Lie algebra action of $u\in U\cong \lie{so}(3)$. But $f_u$ is the unique such function.
\end{prf}

\begin{rmk} Note that over the complex numbers this result still holds. The action of $u$ on the irreducible representation $H$ of $SO(3)$   has a unique one-dimensional kernel whether $u$ is semi-simple or nilpotent.
\end{rmk}

\begin{prp} \label{ico} Let $\pm a_1,\dots, \pm a_6$ be  the vertices of a regular icosahedron, and $X\subset H$ be the space spanned by the harmonic cubics $f_{a_i}$. Then $\omega_u$ restricted to $X$ vanishes for all $u$.

\end{prp}
\begin{prf}
\noindent 
We need to show that
$$\omega_u(f_{a},f_{b})=0$$
for $a,b\in \{ a_1,\dots a_6\}$ or equivalently that
$$(u\cdot f_a,f_b)=(u\cdot f_a)(b)=0.$$
But $u\cdot f_a$ is given by  (\ref{fund})  and, as we saw in Section \ref{clebsch}, $(u\cdot f_a)(b)=0$  for $a=a_{i},b=a_{j}$.
\end{prf}

We call a subspace on which $\omega_u$ vanishes {\it isotropic}.
Over the real numbers, we have a converse to Proposition \ref{ico}.

\begin{prp} \label{icocon} Let $X\subset H$ be a  $3$-dimensional subspace  which is isotropic for all $\omega_u$. Then $X$ is spanned by $f_{a_i}$ where the $a_i$ are the vertices of a regular icosahedron. 
\end{prp} 

\begin{prf} Over $\C$ this is the approach of S.Mukai \cite{Muk}, who described a family of Fano threefolds in terms of  the subvariety  of the Grassmannian $G(3,7)$ of three-dimensional subspaces $X\subset \C^7$ on which a three-dimensional space of skew forms vanishes. Over $G(3,7)$ we have the universal rank $3$ bundle $E$, and each skew form defines a section of $\Lambda^2E^*$, also a rank three bundle. Thus $\omega_1,\omega_2,\omega_3$ define three sections of $\Lambda^2E^*$ and  their vanishing gives a subvariety $F$ of codimension $9$. Since $\dim G(3,7)=12$ this  is three-dimensional. It contains, from Proposition \ref{ico}, the three-dimensional space $SO(3,\C)/A_5$ of (complex) icosahedra and Mukai shows that $F$ is a smooth equivariant compactification of this quotient.  The real points in this open orbit constitute the compact space $SO(3)/A_5$, the Poincar\'e sphere, or space of icosahedra. With these facts, the proof consists of checking that the complement of $SO(3,\C)/A_5$ has no real points. 

The compactification is achieved by adjoining a union of lower-dimensional orbits of $SO(3,\C)$ and in particular each such point is fixed by a one-parameter subgroup.  If $X$ is real this means it is preserved by the Lie algebra action of a vector $u\in U$. But then it is generated by $f_u$, which is  annihilated by $u$, and weight vectors $p$,$\bar p$ where $u\cdot p=\lambda p$, $\lambda\ne 0$. But then 
$$\omega_u(p,\bar p)=(u\cdot p,\bar p)=\lambda(p,\bar p)\ne 0$$
and $\omega_u$ does not vanish on $X$. We deduce that the real subspaces on which $\omega_u$ vanishes are all generated by icosahedra.
  \end{prf}

\section{ Vector bundles}\label{pf}

From the previous section, we see that to find the icosahedra whose vertices lie in the zero set of a harmonic cubic $p$, we must determine  subspaces $X\subset H$  on which the $\omega_u$ vanish. 
Given $p\in H$,  let $W\subset H$ be its $6$-dimensional orthogonal complement. We again work first over $\C$ and then specialize to $\R$. 
\begin{prp} \label{Wdegen} The form $\omega_u$ restricted to $W$ is degenerate if and only if $p(u)=0$.
\end{prp}
\begin{prf} We showed in Proposition \ref{degen} that $f_u$ spans the degeneracy subspace, so $\omega_u$ restricted to $W$ is degenerate if and only if $f_u\in W$, i.e. 
$(p,f_u)=p(u)=0.$
\end{prf}

The skew form $\omega_u$ is linear in $u\in \C^3$ and thus defines a map of bundles on the projective space $\PP(U)$:
$$\omega: W(-2)\rightarrow W^*(-1).$$
The determinant of this is homogeneous of degree $6$, but because $\omega_u$ is skew, it is  the square of a cubic polynomial -- the Pfaffian ($\omega_u^3\in \Lambda^6W^*(3)$). Thus  Proposition  \ref{Wdegen} tells us that  the cubic $p$ is given by the equation of  the Pfaffian of $\omega_u$ on $W$.

\begin{rmk} This derivation of the Pfaffian can be applied to {\it any} irreducible representation space $H$ of $SO(3)$, for, up to a scalar multiple,  there is a unique vector annihilated by an $SO(2)\subset SO(3)$. Given $v\in H$ the Pfaffian defines a homogeneous polynomial of degree $\ell=(\dim H -1)/2$, and so it provides a canonical identification of the projective space of $H$ with the projective space of spherical harmonics of degree $\ell$ .
\end{rmk}

Describing a plane curve as a determinant is a classical problem and this, together with the related Pfaffian problem, is discussed from a modern point of view  in  Beauville's paper  \cite{Bea}. For $u\in C$, 
the map $\omega_u: W(-2)\rightarrow W^*(-1)$ is singular and has a non-zero cokernel. When the curve $C$ defined by $p(u)=0$ is smooth, this defines  a rank $2$ holomorphic vector bundle $E$ over $C$. It has the properties (see  \cite{Bea}) $\det E\cong K_C$ and $H^0(C,E)=0$.

Conversely,  given  any   rank $2$ bundle $E$ on a nonsingular plane curve $C\subset \PP^2$ of degree $d$ with $\det E\cong K_C$ and $H^0(C,E)=0$, then $C$ can be obtained as a Pfaffian, i.e. $E$ admits a natural resolution

$$0\rightarrow {\mathcal O}_{\PP^2}(\C^{2d}(-2))\stackrel{A}\rightarrow {\mathcal O}_{\PP^2}(\C^{2d}(-1))\rightarrow {\mathcal O}_C(E)\rightarrow 0$$

where $A$ is skew-symmetric and linear in $u$ and $C$ is defined by the vanishing of the Pfaffian of $A$. 

\vskip .25cm
In our case, $d=3$, and $K_C$ is trivial so we have a rank $2$ vector bundle $E$ with $\det E$ trivial and $H^0(C,E)=0$. As Atiyah showed in \cite{At}, when $C$ is smooth, such a bundle is either of the form
$$E=L\oplus L^*$$ 
for a non-trivial degree zero  line bundle $L$, or is a non-trivial extension $$0\rightarrow L\rightarrow E\rightarrow L\rightarrow 0$$ of a non-trivial line bundle $L$ with $L^2\cong{\mathcal O}$. 
\vskip .25cm
\begin{prp} \label{vb} For a smooth cubic curve $C$, there is a one-to-one correspondence between $3$-dimensional subspaces $X\subset W$ on which $\omega_u$ vanishes for all $u$ and degree $0$ holomorphic line bundles $L\subset E^*$.
\end{prp}
\begin{prf} Suppose that  $\omega_u$ vanishes on $X\subset W$.  The skew form $\omega_u$ defines a homomorphism 
$$\omega: W\rightarrow W^*(1)$$
whose kernel is $E^*(-1)$.  Restricting to $X$, we have
$$\omega\vert_{X}: X\rightarrow X^o(1)$$
where $X^o\in W^*$ is the annihilator of $X$. 

The cubic $p$ is defined by $\det \omega\vert_{X}=0$, and so 
as $u$ varies over $C$, the kernel of $\omega\vert_{X}$ describes a line bundle $L(-1)\subset X$ where $\deg L=0$. Since this is also in the kernel of $\omega:  W\rightarrow W^*(1)$ it follows that  $L(-1)\subset E^*(-1)$. 
\vskip .25cm
Conversely, suppose $L(-1)\subset E^*(-1)\subset W$. From the exact sequence of bundles on $C$ 
$$0\rightarrow E^*(-2)\rightarrow W(-1)\rightarrow W^*\rightarrow E(1)\rightarrow 0$$
we have $W^*\cong H^0(C,E(1))$. Since $L^*(1)$ has degree $3$,  $H^0(C,L^*(1))\subset H^0(C,E(1))$ defines a three-dimensional subspace $X^o\subset W^*$ and $X\subset W$. Moreover, since a line bundle of degree three embeds a smooth elliptic curve in the plane, we have $L(-1)\subset X$.

The restriction of $\omega_u$ to $X$ is a section $w$ of $\Lambda^2X^*(1)$. Identifying $\Lambda^2X^*\cong X\otimes \Lambda^3 X^*$, we see that this defines a section $w^*$ of $X(1)$ which at each point lies in the degeneracy subspace of the skew form $w$ on $X$. But $L(-1)\subset E^*(-1)$ and so $L(-1)$ is the degeneracy subspace of $\omega_u$ on $X$. Thus if $w^*$ is non-zero it defines a section of $L$. But in the Pfaffian construction $E$, and hence $L$,  has no holomorphic sections, so we have a contradiction unless $w^*=0$, i.e. $\omega_u$ vanishes on $X$.
\end{prf}

Now when $E\cong L\oplus L^*$ and $L^2$ is non-trivial, $L$ and $L^*$ are the only degree zero  subbundles. For the nontrivial extension $0\rightarrow L\rightarrow E\rightarrow L\rightarrow 0$, $L$ is the only subbundle. The bundle $E=L\oplus L$ for $L^2$ trivial clearly has infinitely many subbundles isomorphic to $L$. Hence, 
using Propositions \ref{icocon} and \ref{vb}  we see that for a nonsingular cubic  Atiyah's theorem offers the alternatives:

\begin{itemize}
\item
If $E\cong L\oplus L^*$ there are two isotropic subspaces $X_1,X_2\subset W$. 
If both $X_1$ and $X_2$ are real,  $p$ vanishes on the vertices of two  icosahedra. If  $X_1=\bar X_2$  there are no real icosahedra.  
\item
When $E$ is a non-trivial extension   with $L$ non-trivial and $L^2\cong{\mathcal O}$, then there is a unique icosahedron.
\item
When $E=L\oplus L$ there are infinitely many icosahedra.
\end{itemize}

The discussion above is for smooth curves, but the situation for nodal cubics or reducible ones with transverse intersections is  similar \cite{Burb}.  Cuspidal cubics present far more problems, but fortunately we do not have to deal with these. In fact, as noted in \cite{Ley} and proved in \cite{Cheng}, singularities of nodal sets are quite simple. If $\Delta f= \lambda f$ and $f$ together with its derivative vanish at a point, then the function is locally approximated by a solution to $\Delta f=0$, i.e. the real part of $cz^n$. So the  nodal set at a singularity locally consists of $n$ smooth curves meeting at the same angle $\pi/n$: for example the three great circles given by $x_1x_2x_3=0$  meet at right angles. 

Another regularity issue is that for a value of $u$ at which $\omega_u$ restricted to  $W$ is degenerate, it has rank $4$ and no lower. This is because $(u\cdot f,q)=0$ for all $q\in W$ if and only if $u\cdot f=\lambda p$. But if 
$u\cdot f_i=\lambda_i p$ for $i=1,2$ then $u\cdot (\lambda_2f_1-\lambda_1 f_2)=0$ so $\lambda_2f_1-\lambda_1 f_2=\mu f_u$. Hence there is only a two-dimensional space of polynomials $f$ with $(u\cdot f,q)=0$ for all $q\in W$. In particular this means that even when the curve $C$ is singular the degeneracy subspace still gives  a well-defined vector bundle. 

On a nodal curve, one can understand vector bundles by passing to the normalization -- a copy of $\PP^1$  with two distinguished points  which map to the node. If the pull-back of $E$ is trivial, then the identification of the fibres at the two points  is given by  a matrix  $A\in SL(2,\C)$ and if $H^0(C,E)=0$, $A$ must not have $1$ as an eigenvalue. The conjugacy classes are then
$$\pmatrix{a & 0\cr
             0 & a^{-1}} a \ne \pm 1,\quad \pmatrix{-1 & b\cr
             0 & -1} b\ne 0,\quad  \pmatrix{-1 & 0\cr
             0 & -1}$$
             corresponding exactly to Atiyah's classification.
If the bundle is non-trivial on $\PP^1$, since $H^0(C,E)=0$, it can only be isomorphic to ${\mathcal O}(1)\oplus {\mathcal O}(-1)$. In this case there is no degree zero subbundle, so from Proposition \ref{vb}, there is no icosahedron, real or complex, whose axes lie  on $C$. 

We could proceed in this manner to the reducible cubics but these can be dealt with more concretely as we shall see in the next section.

\section{Counting icosahedra}

Given a harmonic cubic, we  want to determine how many icosahedra have all their vertices on its nodal set.  For a nonsingular or nodal cubic the previous section gives us the information that this number  is zero, one, two or infinitely many. There is an elementary way of going the other way round -- associating a harmonic cubic to a pair of icosahedra:

\begin{prp} \label{two} Given two  regular icosahedra with vertices on the unit sphere, then 
\begin{itemize}
\item
if they have no vertices in common, there exists a one-dimensional space of   harmonic cubics  which vanish at all vertices,
\item
 if they have one vertex $a$ in common then there is a two-dimensional space of such  cubics, each  of the form  
$$(u,x)((a,x)^2-\frac{1}{5}(x,x)(a,a))$$
where $(u,a)=0$.
\end{itemize}
\end{prp}
\begin{prf} The two icosahedra define two $4$-dimensional subspaces $V_1$ and $V_2$ in the $7$-dimensional space $H$, which must therefore intersect non-trivially. 

Suppose that $\dim V_1\cap V_2>1$, then we have a pencil of plane cubic curves passing through the  points $[a_1],\dots, [a_6]$, $[b_1],\dots,[b_6]$ defined by the axes of  the two icosahedra. If two axes coincide then so do the icosahedra, so at most one axis can be common. This means the curves intersect in at least eleven distinct points. But by B\'ezout's theorem, there can be a maximum of nine  unless there is a common component. 

Suppose that component is a line. No three of the $[a_{i}]$ are collinear (or the $[b_{i}]$) so there is a maximum of four ({\it e.g.} $[a_1],[a_2],[b_1],[b_2]$) on a line, leaving at least seven. The pencil of conics which remains has a maximum of four intersections, so that must have a common line. Again that can have at most four new vertices, leaving  three for a pencil of lines, which is a contradiction.

If the common component is a conic, then it contains at most five axes from either icosahedron (as we saw in the proof of Proposition \ref{lem}.) The remaining axis must lie in the pencil of lines. If $[a]$ is that axis, then the other five axes make the same angle with $a$ and so satisfy the equation $$(a,x)^2-\frac{1}{5}(x,x)(a,a)=0.$$ This is the fixed conic. The remaining line has equation $(u,x)=0$ and since it passes through $[a]$, $(u,a)=0$. 
\end{prf}

We see here that a pair of icosahedra without a common axis determines, up to a scalar multiple, a unique harmonic cubic. We show next:

\begin{prp} \label{three} Suppose $p$ vanishes on the vertices of more than two icosahedra. Then $p$  is of the form 
$$(u,x)((a,x)^2-\frac{1}{5}(x,x)(a,a))$$
with $(u,a)=0$. From  Proposition \ref{two} it vanishes for infinitely many icosahedra, all with vertex $a$.
\end{prp}

\begin{prf} Let $X_1,X_2$ be  isotropic subspaces of $W$ for $\omega_u$ and assume that $p$ is not of the form in Proposition \ref{two}. Then from that Proposition, $\dim V_1\cap V_2=1$ and hence $X_1\cap X_2=0$. A third subspace $X_3$ intersects each of these trivially and hence is the graph of an invertible  linear transformation $S:X_1\rightarrow X_2$. Thus for $w_1,w_2\in X_1$ we have
$$0=\omega_u(w_1+Sw_1,w_2+Sw_2)=\omega_u(w_1,Sw_2)+\omega_u(Sw_1,w_2)$$
since $\omega_u(w_1,w_2)=0=\omega_u(Sw_1,Sw_2)$. This means the graph of $tS$ for any real number $t$ is also  isotropic, so we have a one-parameter family of icosahedra.

Since  icosahedra are all in the same $SO(3)$ orbit, differentiating at $t=0$ we see that $Sw=u\cdot w$ for some $u$. Moreover, since $S$ is invertible, $u$ is not a multiple of the vertices $a_i$ of the icosahedron defined by $X_1$. But this means that each element of $X_2=u\cdot X_1$ is orthogonal to $X_1$ because 
$$(u\cdot w_1,w_2)=\omega_u(w_1,w_2)=0.$$
In particular if $b$ is a vertex of the icosahedron for $X_2$, $f_b\in X_2$ is orthogonal to the $f_{a_i}\in X_1$, or equivalently $f_b(a_i)=0$. 

However, as remarked in the introduction, no $f_u$ can vanish on the vertices of an icosahedron, so we have a contradiction.

\end{prf}

It remains to discuss the other  reducible cubics which pass through the axes of the icosahedron. First consider the  case of a line through two axes and a nonsingular conic through the remaining four. In the standard model if we take the line to be $x_3=0$, then the conic is of the form 
\begin{equation}
a(\phi^4x_1^2+x_2^2-\phi^2x_3^2)+bx_1x_2=0
\label{conic1}
\end{equation}
where $x_1^2+x_2^2+x_3^2=1+\phi^2$ and for nonsingularity $a\ne 0, b^2\ne 4a^2\phi^4$. 

Now rotate the icosahedron an angle $\theta$ about the $x_3$-axis. Four vertices still lie in $x_3=0$, and we want to see if the remaining eight lie on the conic. By symmetry it is enough to consider 
$$(\phi \sin \theta, \phi \cos\theta,1)\qquad (\cos\theta,-\sin\theta,\phi).$$ 
Substituting in (\ref{conic1}) this gives either $\sin \theta=0$ -- the original icosahedron, or
$$\tan \theta=\frac{b}{a(1-\phi^4)}.$$
There are thus  two icosahedra (and no more by  Proposition \ref{three}).

Now consider three lines passing through the six points. There are $15$ possibilities which break up into  two $A_5$ orbits. The first  are the five cubics $q_1,\dots,q_5$ we started with, and we observed in the introduction that they vanish on the vertices of two icosahedra (this is also the case $a=0$ in (\ref{conic1})). The other $10$ are of the form  $q_{\alpha}-q_{\beta}$ where the cubic curve consists of three lines meeting at  a point. On the sphere these are three equally spaced lines of longitude. Since no three axes in $\PP(U)^r$ are collinear, we must have $4$ vertices in each great circle and the North-South axis passes through the midpoint of a face of the icosahedron. But this determines the distances of the vertices from the  poles, so the icosahedron is unique.

In discussing these special cases, we have shown that on the line in $\PP(V)$ defined by the family (\ref{conic1}) there are three points $a=0, b=\pm 2a\phi$  where there is a unique icosahedron:  the three degenerate cubics consisting of $L_{12}$  together with $\{L_{34},L_{56}\},\{L_{35},L_{46}\}$ or $\{L_{36},L_{45}\}$. We shall see this more generally next. 

\section{The Clebsch cubic revisited}

 \begin{prp} A real harmonic cubic $p\in V$ vanishes on the vertices of precisely two distinct icosahedra if and only $[p]$ does not lie on the Clebsch cubic surface. 
\end{prp}
\begin{prf} Let $X\subset H$ be the $3$-dimensional subspace defining an  icosahedron and $p$ be a cubic which vanishes at its vertices. Consider as in Proposition \ref{vb} $\omega\vert_{X}: X\rightarrow X^o(1)$, whose determinant vanishes on the plane cubic $C\subset \PP(U)$ defined by $p$.  The kernel $L(-1)$ lies in the kernel $E^*(-1)$ of $\omega: W\rightarrow W^*(1)$. From Atiyah's classification, for a smooth or nodal cubic there is a unique icosahedron if and only if  $L^2$ is trivial but $L$ is non-trivial. 

As $u$ varies in $C$ consider $f_u$. Now $\omega_u(f_u,q)=(u\cdot f_u,q)=0$ so that $f_u$ generates a line bundle in the kernel of $\omega: W\rightarrow W^*(1)$. We identify this bundle as follows. Consider the map
$f:\PP(U)\rightarrow \PP(H)$ defined by $u\mapsto f_u$. Then $f_u$ spans the pull-back $f^*({\mathcal O}(-1))$ and since $f_u$  (see (\ref{harm})) is homogeneous of degree $3$ in $u$, $f^*({\mathcal O}(-1))\cong{\mathcal O}(-3)$. Restricted to $C$, this is the line bundle we want. 

We thus have an inclusion 
${\mathcal O}(-3)\subset E^*(-1)$ over $C$. Projecting to  $E^*(-1)/L(-1)$ gives a section of $L^*(2)$. This vanishes when $f_u\in X$ which from Proposition \ref{ico}  is at the six points $[a_1], \dots, [a_6]$ defined by the axes of the icosahedron. Thus, on $C$ we have the relation of divisor classes 
\begin{equation}
L\sim {\mathcal O}(2)- \sum_1^6[a_i].
\label{divclass}
\end{equation}

Now take three general points on $C$ distinct from the $a_i$ and blow up $\PP(U)$ at these $9$ points. We get an elliptic surface $Y$: the anticanonical bundle $K_Y^*$ has two sections which define a projection $\pi:Y\rightarrow \PP^1$ with elliptic fibres, which consist of the pencil of cubics through the nine points.  Let $E_i$, $1\le i\le 9$ be the exceptional curves on $Y$, then  

\begin{equation}
K_Y\sim -3 H+\sum_1^9 E_i
\label{antican}
\end{equation}

where $H$ is the pull-back of the hyperplane divisor on $\PP(U)$. 
Let $J$ be the line bundle on $Y$ whose divisor class is 
$$J\sim  2H-\sum_1^6 E_i$$
then from (\ref{divclass}) $J$ restricts to $L$ on each elliptic fibre of $\pi$. 

\begin{rmk} This formula identifies $L$, which has come to us from the consideration of vector bundles, with a more basic geometric object. The elliptic curve $p(x)=0$ lies in two projective planes: $\PP(U)$ and (via its proper transform) a plane section of $S\subset \PP(V)$. The two different hyperplane bundles differ by $L$.
\end{rmk}

We shall apply Grothendieck-Riemann-Roch to $J^2$ on $Y$ to determine the number of cubics in the pencil for which $L^2$ is trivial.

 Now $\pi_*[K_S]=0$ since it is dual to a fibre, and for the exceptional curves  $\pi_*[E_i]=1.$   Hence
$\pi_*[H]=3$ and $p_*[J]=0$. Moreover we have intersection numbers
$[J]^2=-2, [K_Y][J]=0$. Hence
$$\pi_*(\td(Y) \ch(J^2))=-3[\PP^1].$$
Hence $\ch(\pi_{!}(J^2))=-3[\PP^1]$ 
 and so $L^2$ is trivial  for three cubics in the pencil. 
 
 It follows that a generic line in $\PP(V)$ meets the locus of cubics for which $L^2$ is trivial in three points, so that locus is a cubic surface.  But it is invariant under the icosahedral group and hence must be the Clebsch cubic $S$. 

All real cubics in $\PP(V)$ pass through the axes of the standard icosahedron, moreover from Proposition \ref{three} those that pass through infinitely many form  six of the $27$ lines on $S$. We also treated the reducible cases in the previous section, so any  $[p]$ in the complement of $S$ must pass through the axes of two distinct icosahedra and no more.

 \end{prf}
 \section{The general case}
 So far we have mainly focused on the four-dimensional  space $V$ of harmonic cubics which vanish on the vertices of a fixed icosahedron, and identified those which vanish on a second. We now ask the main question of the seven-dimensional space $H$ of all harmonic cubics: which ones vanish on the vertices of {\it two}  icosahedra? In what follows we are only considering real points, so we omit the superscript $r$ to denote this. 

Clearly if $p$ vanishes for icosahedra $I_1,I_2$, we can transform $I_1$ to the standard icosahedron $I$ by the action of $SO(3)$, so that the polynomials we are seeking are the transforms of the complement of $S$ under the action map
$$F: SO(3)\times \PP(V)\rightarrow \PP(H).$$
The points of $\PP(H)$ for which there is either one or infinitely many icosahedra are given by the orbit of the cubic surface $S$, which is an $SO(3)$-invariant  hypersurface in $\PP(H)$. Now  
 the ring of $SO(3)$-invariants of the $7$-dimensional representation $H$  is  generated in degrees
$2, 4, 6, 10,$ and $15$, with the square of the one of degree $15$ being a polynomial
in the other four generators, which are themselves algebraically
independent. (I thank Robert Bryant for this information). There is  no  invariant of degree three but in the Appendix we calculate a sextic invariant $J$ on $H$  which, when restricted to $V$, becomes $\sigma_3^2$. Consequently for real harmonic cubics, the orbit of $\PP(V)$ in $\PP(H)$ under the action of   $SO(3)$ is contained in the subset $J(f)\ge 0$.

 \begin{prp} If $J(f)>0$ then $[f]$ is in the orbit of $\PP(V)\setminus S$ under $SO(3)$.
 \end{prp}
 \begin{prf} From Proposition \ref{two}, the orbit of $\PP(V)\setminus S$ under the $SO(3)$ action is the set of unordered pairs of icosahedra having no axis in common. Similarly  $\PP(V)\setminus S$ consists of icosahedra having no axis in common with the standard one $I$. The ordering of the two icosahedra is the basis for what follows.

Define an equivalence relation on $SO(3)\times \PP(V)$ by $(g,p)\sim (h,q)$ if $[gp]=[hq]$ and $gI=hI$. Let $M$ be the space of equivalence classes.  Clearly $F(g,p)=[gp]$ factors through $M$. Since $h^{-1}g\in A_5$ the equivalence class of $(g,p)$ is finite and in one-to-one correspondence with the subgroup of $A_5$ fixing $[p]$. We may see these fixed points directly by looking at the permutation representation and in all cases the fixed point set in $\PP(V)$ is a line, and so of codimension $2$. So $M$ is a six-manifold.

If $[p]\notin S$ there are two icosahedra and so given $[gp]$, two equivalence classes. This represents $F$ as a double covering, and the involution interchanging the sheets extends  to $M$ with fixed point set the five-dimensional image of $SO(3)\times S$. Thus  $F$ defines a map from an (unoriented) six-manifold $\bar M=M/\Z_2$ with boundary  to the six-dimensional $\PP(H)$ where the boundary maps to the hypersurface $Y$ defined by $J(f)=0$. The image of $F$ lies in the region $\PP(H)^+$ where $J\ge 0$. Since $F$ is one-to-one on the interior, it is of degree one and hence the map 
$$F_*:H_6(\bar M,\partial \bar M;\Z_2)\rightarrow  H_6(\PP(H)^+\!\!,Y;\Z_2)$$
is an isomorphism. In particular it follows that $F$ is surjective, and so $J(f)\ge 0$ means that $f=gp$ for $[p]\in \PP(V)\setminus S$. Hence  $f$ vanishes on the vertices of  exactly two icosahedra. 
\end{prf} 
We thus obtain the theorem:

\begin{thm}  Let $f$ be a degree three  spherical harmonic.
\begin{itemize}
\item  If $J(f)>0$ (resp.  $J(f)<0$) then the nodal set of $f$ contains the vertices of two (resp. zero) regular icosahedra. 
\item
When $J(f)=0$ and $f$ is of the form  $f(x)=(u,x)((a,x)^2-(x,x)(a,a)/5)$
with $(u,a)=0$ then any regular icosahedron with $a$ as a vertex lies on the nodal set, otherwise there is a unique such icosahedron.
\end{itemize}
\end{thm}

\section{Appendix: Invariants}\label{invariant}

We calculate explicitly here the $SO(3)$-invariant polynomial $J$ on $H$ which restricts to $\sigma_3^2$ on $V$.
\vskip .25cm

Let $f(x)$ be a spherical harmonic of degree $3$ and consider the symmetric form $S$ on $U$ defined by 
$$S(u,v)=\int_{S^2}\!\!\!f^2(x)(u,x)(v,x)\,dS$$
or
$$S_{ij}=\int_{S^2}\!\!\!f^2(x)x_ix_j\,dS$$
The $SO(3)$-invariants $\tr S,\tr S^2,\tr S^3$ define invariants of degree $2,4,6$ of $f$, and we shall relate these to the $A_5$-invariants when we restrict to our distinguished four-dimensional space of  harmonic cubics of the form 
$$f(x)=\sum_1^5y_{\alpha}q_{\alpha}(x)$$
where $\sum_{\alpha} y_{\alpha}= 0$ and $\sum_{\alpha} q_{\alpha}=0$.

First consider $f=q_{\alpha}$. In this case $S$ is invariant under conjugation by the stabilizer $A_4\subset A_5$ of $\alpha$. But $A_4$, the tetrahedral group, fixes no axes in $\R^3$, so that 
\begin{equation}
\int_{S^2}\!\!\!q_{\alpha}^2x_ix_j\,dS =\lambda \delta_{ij}
\label{X1X1}
\end{equation}

Taking the trace
$$3\lambda= \int_{S^2}\!\!\!q_{\alpha}^2\,dS.$$
To evaluate this and similar integrals, note that  each $q_{\alpha}$ is of the form $(a_1,x)(a_2,x)(a_3,x)$, and by the invariant theory of the orthogonal group
$$\int_{S^2}(a_1,x)(a_2,x)\dots(a_{2n},x)\,dS=c_n\sum(a_{i_1},a_{i_2})\dots (a_{i_{2n-1}},a_{i_{2n}})$$
for a universal constant $c_n$ where the sum is over partitions $\{i_1,i_2\}\dots\{i_{2n-1},i_{2n}\}$ of $\{1,2,3,\dots,2n\}$ into two-element subsets. 
In fact, taking all the $a_i=(0,0,1)$ we  evaluate   $c_n=2^{n+2}\pi n!/(2n+1)!$

If $q_{\alpha}(x)=x_1x_2x_3=(e_1,x)(e_2,x)(e_3,x)$ then by orthogonality there is only one partition of $\{e_1,e_2,e_3,e_1,e_2,e_3\}$ which gives a non-zero result. Thus 
\begin{equation}
\lambda=\frac{1}{3}c_3
\label{lam}
\end{equation} 

Now consider 
$$\int_{S^2}\!\!\!q_{\alpha}q_{\beta}x_ix_j\,dS.$$
This is invariant by the stabilizer $S_3$ of the unordered pair $\{\alpha,\beta\}$.

The (dihedral) action of $S_3$ has one invariant axis, which we represent by a unit vector $v^{\alpha\beta}$ and then 
\begin{equation}
\int_{S^2}\!\!\!q_{\alpha}q_{\beta}x_ix_j\,dS=\mu \delta_{ij}+\nu v^{\alpha\beta}_iv^{\alpha\beta}_j
\label{X1X2}
\end{equation}
\begin{lemma} $\mu=0$, $\nu=-3\lambda/4$.
\end{lemma}
\begin{lemprf}
Taking the trace of (\ref{X1X2}) gives 
$$3\mu+\nu=\int_{S^2}\!\!\!q_{\alpha}q_{\beta}\,dS$$
and using $\sum_{\alpha} q_{\alpha}=0$ we obtain
$$\int_{S^2}\!\!\!q_{1}q_{5}\,dS=-\int_{S^2}\!\!\!q_{1}^2+q_{1}q_{2}+q_{1}q_{3}+q_{1}q_{4}\,dS$$ and so
\begin{equation}
3\lambda+12\mu+4\nu=0.
\label{lmn}
\end{equation}

The action of $A_5$ on the $10$ unordered pairs $\{\alpha,\beta\}$ is, in the icosahedron, the action on opposite pairs of the $20$ faces, and $v^{\alpha\beta}$ is the axis joining their centres. In Figure 1 such an  axis points towards the reader, and passes through the centre of the triangle $ABC$.  In Figure 2 it defines an Eckard point.

The three shaded  planes in the diagram form the zero locus of $q_{\alpha}$. The  face $ABC$  has three adjacent faces and  each shaded plane provides one of its edges. The remaining three edges define the second triple of planes which is the zero locus of $q_{\beta}$.

Analytically suppose that $q_{\alpha}(x)=x_1x_2x_3$ and $v=v^{\alpha\beta}=(1,1,1)/\sqrt{3}$. The vertices $A,B,C$  are then $(0,\phi,1)$, $(\phi,1,0)$, $(1,0,\phi)$. Then $q_{\alpha}$ and $q_{\beta}$ are interchanged by a rotation by $\pi$ in $S_3$ which takes $v$ to $-v$ and the vertex $(0,\phi,1)$ to its opposite $-(0,\phi,1)$. This is the rotation
\begin{equation}
g=\frac{1}{2}\pmatrix {-\phi & \phi-1&-1\cr
                         \phi-1 & -1 & -\phi\cr
                         -1 & -\phi & \phi-1}
                         \label{gmatrix}
                         \end{equation}

Now from  (\ref{X1X2}) we see that 
$$S(v^{\alpha\beta},v^{\alpha\beta})=\int_{S^2}\!\!\!q_{\alpha}q_{\beta}(v^{\alpha\beta},x)^2\,dS=\mu +\nu$$
and we can evaluate this  by partitioning
$\{e_1,e_2,e_3,ge_1,ge_2,ge_3,v,v\}$ into pairs.

If $v$ is paired with $v$ then since $(v,v)=1$ we get as the sum of these contributions 
\begin{equation}
c_3^{-1}\int_{S^2}\!\!\!q_{\alpha}q_{\beta}\,dS=c_3^{-1}(3\mu+\nu)=-\frac{1}{4}\label{cxx}
\end{equation}
from (\ref{lmn}) and (\ref{lam}).
 
If both $v$ terms are paired with $e_i$ terms or $ge_j$ terms then by orthogonality one of the remaining inner products is zero. Thus the only non-zero contributions to the integral  have factors
$(v,e_i)(v,ge_j)$ which (since $gv=-v$) is $-(v,e_i)^2=-1/3$. This gets multiplied by terms $(e_k,ge_l)(e_m,ge_n)$. These are products of entries in  (\ref{gmatrix}) and the sum of these $18$ terms is readily evaluated to be $3$. Because of the repeated terms $v,v$ there are two partitions which give the same contribution, and so we find these terms producing a  contribution
$-2$. Adding in (\ref{cxx}) we find
$$\mu+\nu=-\frac{9}{4}c_4.$$
But $3\mu +\nu=-3\lambda/4=-c_3/4$ and $c_3/c_4=9$, so $3\mu +\nu=\mu+\nu$ and therefore $\mu=0$ and $\nu=-3\lambda/4$. 
\end{lemprf}

From the lemma we may as well normalize $S$  and take $\lambda=4,\mu=0,\nu=-3$.
Consider a general $f=\sum_{\alpha}y_{\alpha}q_{\alpha}$. Then 
\begin{eqnarray*}
S_{ij}&=&\sum_{\alpha,\beta}\int_{S^2}\!\!\!y_{\alpha}y_{\beta}q_{\alpha}q_{\beta}x_ix_j\,dS\\
&=&\lambda \sum_1^5y_{\alpha}^2 \delta_{ij}+\sum_{\alpha\ne\beta}y_{\alpha}y_{\beta}(\mu \delta_{ij}+\nu v^{\alpha\beta}_iv^{\alpha\beta}_j)\\
&=&-8\sigma_2\delta_{ij}-3\sum_{\alpha\ne\beta}y_{\alpha}y_{\beta} v^{\alpha\beta}_iv^{\alpha\beta}_j
\end{eqnarray*}
using the elementary symmetric polynomials $\sigma_i$ in $y_{\alpha}$ and   the condition $\sigma_1=0$.
In particular $\tr S=-30\sigma_2$.
 Define the symmetric matrix $M$ by 
$$M_{ij}=\sum_{\alpha\ne\beta}y_{\alpha}y_{\beta} v^{\alpha\beta}_iv^{\alpha\beta}_j$$
then we have 
\begin{eqnarray*}
\tr M^2&=&\sum_{\alpha\ne\beta}y_{\alpha}y_{\beta} \sum_{\gamma\ne\delta}y_{\gamma}y_{\delta} (v^{\alpha\beta},v^{\gamma\delta})^2 \\
 \tr M^3 &= &\sum_{\alpha\ne\beta}y_{\alpha}y_{\beta} \sum_{\gamma\ne\delta}y_{\gamma}y_{\delta} \sum_{\mu\ne\nu}y_{\mu}y_{\nu}(v^{\mu\nu},v^{\alpha\beta})(v^{\alpha\beta},v^{\gamma\delta})(v^{\gamma\delta},v^{\mu\nu})
\end{eqnarray*}
and for this we need to know the inner product terms involving the $v^{\alpha\beta}$.  

As we noted, the action of $A_5$ on  unordered pairs $\{\alpha,\beta\}$ is the action on opposite pairs of faces. There are $15$ pairs of these $\{\{\alpha,\beta\},\{\gamma,\delta\}\}$ that have no number in common, and these correspond to  {\it adjacent} faces -- or, considering the common edges, to opposite pairs of the $30$ edges.  The other $30$ pairs $\{\{\alpha\beta\},\{\beta\gamma\}\}$ correspond to pairs with no common edge.  In Figure 1 $D=(0,\phi,-1)$, so a unit vector passing through the centre of triangle $ADB$ is $(\phi,2\phi,0)/\sqrt 3 (\phi+1)$. Since $ABC, ADB$ are adjacent this  gives
$(v^{\alpha\beta},v^{\gamma\delta})^2= 5/9$. The vertex $E$ is $(1,0,-\phi)$ and $BDE$ is not adjacent to $ABC$ or its opposite. The  vector through the centre of  triangle $BDE$ is  $(1,1,-1)/\sqrt{3}$ so here we obtain $(v^{\alpha\beta},v^{\beta\gamma})^2= 1/9$. 

We can now evaluate 
$\tr M^2$ in terms of elementary symmetric functions $\sigma_i$. Using $\sigma_1=0$ we obtain
$$\tr M^2=4\sigma_2^2+\frac{160}{9}\sigma_4.$$

To evaluate $\tr M^3$ we need to consider triples of two-element subsets of $\{1,\dots,5\}$. There are $5$ orbits represented by 
$$\{1,2\}\{1,2\}\{1,2\}\quad \{1,2\}\{1,2\}\{1,3\}\quad \{1,2\}\{1,2\}\{3,4\}\quad\{1,2\}\{1,3\}\{1,4\}$$
and $\{1,2\}\{1,3\}\{4,5\}$
where the last one  contributes zero since it occurs with a coefficient given by the symmetric polynomial $\sigma_1\sigma_5$ and in our case $\sigma_1=0$. 

The orbit of $\{1,2\}\{1,3\}\{1,4\}$ introduces a new expression $(v^{12},v^{13})(v^{13},v^{14})(v^{14},v^{12})$ which we must calculate. To do this, note that $\{1,2\}\{1,3\}\{1,4\}$ is stabilized by the rotation $g$ of order $3$ which fixes $1$ and $5$ and  permutes $2,3$ and $4$. Thus 
$$(v^{12},v^{13})(v^{13},v^{14})(v^{14},v^{12})=(v^{12},gv^{12})(gv^{12},g^2v^{12})(g^2v^{12},v^{12})=(v^{12},gv^{12})^3.$$
If $\{1,5\}$ defines the face $ABC$, then the rotation is $g(x,y,z)=(y,z,x)$ so that taking $v^{12}=(1,1,-1)/\sqrt{3}$ to pass through the centre of $BDE$, we obtain
$(v^{12},gv^{12})^3=-1/27$. 

We now calculate  
$$\tr M=2\sigma_2,\quad  \tr M^2=4(\sigma_2^2+\frac{40}{9}\sigma_4),\quad \tr  M^3=8(\sigma_2^3-\frac{20}{9}\sigma_2\sigma_4+2\sigma_3^2).$$
Denoting by $\tau_i$ the symmetric functions in the eigenvalues of $M$, we obtain

$$16\sigma_3^2=3\tau_3-4\tau_1\tau_2$$

and this gives explicitly the invariant $J$.

\vskip 1cm
 Mathematical Institute, 24-29 St Giles, Oxford OX1 3LB, UK
 
 hitchin@maths.ox.ac.uk

 \end{document}